\input amstex
\documentstyle{amsppt}
\magnification 1200

\topmatter
\title
On convergence to non normal laws for ergodic random fields
\endtitle

\author
Dalibor Voln\'y
\endauthor

\affil
CNRS and University of Rouen Normandy
\endaffil

\abstract
For $\Bbb Z^d$ actions equipped with a completely commuting filtration, normalized partial sums of
martingale differences converge in distribution (see \cite{V19}) but even if the action is ergodic, the
limit law need not be normal. 
Explication of this phenomenon has been studied in \cite{GLV} but often in dimension 2 only. It turned out that 
the non normality is related to the factor of product type $\Cal I$.

Here we extend the necessary and sufficient condition of normality from dimension 2 to all finite dimensions,
as well as a sufficient condition of normality for martingale differences in the orthocomplement of $L^2(\Cal I)$.
At the same time an improved and corrected version of the central limit theorem for random fields from \cite{V19}
with dimension $d>2$ is given.
\endabstract

\endtopmatter

\document

\heading
1. Introduction and main results
\endheading


The probably first CLT for an ergodic sequence of square integrable martingale differences was found by
P\. Billingsley in 1961 (cf\. \cite{B}), independently (by another method and with rate of convergence) by I.A\. Ibragimov in 1963 (\cite{I}).
The results were followed by many limit theorems using martingale approximation, starting with M.I\. Gordin's 
seminal paper \cite{Go}. Many of the results appear in the monography \cite{HaHe}, more recent results can be 
found e.g\. in \cite{MPU}.

The case of random fields has been more complicated, in particular because there are different ways how to 
define the filtration. Here we use {\it completely commuting filtrations} (a definition will be given later). 
While for ergodic sequences of martingale
differences the CLT gives convergence to a normal law, for ergodic random fields this need not be the case. 
The paper tries to understand this phenomenum, continuing in the research presented in \cite{GLV}. 

Let $(\Omega, \Cal A, \mu)$ be a probability space with a $\Bbb Z^d$ action $T_{\underline{i}}$,
$\underline{i} = (i_1,\dots,i_d) \in \Bbb Z^d$. For each $\underline{i}$, $T_{\underline{i}}$ is thus a bijective, 
bimeasurable, and measure preserving mapping of $\Omega$ onto itself, $T_{\underline{i}+\underline{j}} = 
T_{\underline{i}}\circ T_{\underline{j}}$. For $f$ measurable, $f\circ T_{\underline{i}}$ will be sometimes noted
as $U_{\underline{i}} f$. A sub-sigma-algebra $\Cal C$ of $\Cal A$ invariant for all $T_{\underline{i}}$ will be
called a factor; it is a particular case of a factor as defined in ergodic theory. \newline
We say that $(\Cal F_{\underline{i}})_{\underline{i}\in \Bbb Z^d}$ is a {\it completely commuting filtration} of 
sub-sigma-algebras of $\Cal A$ if
\roster
\item For $\underline{i} \prec \underline{j}$ it is $\Cal F_{\underline{i}} \subset \Cal F_{\underline{j}}$
where $\underline{i} \prec \underline{j}$ means $i_1\leq j_1, \dots, i_d\leq j_d$.
\item $\Cal F_{\underline{i}} = T_{-\underline{i}} \Cal F_{\underline{0}}$ where $\underline{0} = (0,\dots,0)$.
\item For $f$ integrable, 
$$
  E\big( E(f | \Cal F_{\underline{i}}) \,\big|\, \Cal F_{\underline{j}} \big) = 
  E\big( E(f | \Cal F_{\underline{i}\wedge \underline{j}})
$$
where $\underline{i}\wedge \underline{j} = (i_1\wedge j_1, \dots, i_d\wedge j_d)$, $u\wedge v = \min(u, v)$.
\endroster

Denote $e_q$ the vector in $(i_1,\dots, i_d) \in \Bbb Z^d$ where all coordinates except $i_q$ are zero, $i_q=1$
($1\leq q\leq d$).
An integrable function $f$ is a {\it martingale difference} if it is a martingale difference for every $T_{e_q}$
and the filtration $(\Cal F_i^{(q)})_i$ where $\Cal F_i^{(q)} = \Cal F_{\underline{j}}$, $\underline{j} =(j_1,\dots,j_d)$,
$j_q=i$ and $j_r=\infty$ for $r\neq q$, $1\leq r\leq d$,
$1\leq q\leq d$. This means, $f$ is $\Cal F_{\underline{0}}$-measurable and for all $1\leq q\leq d$,
$E(f | \Cal F_{-e_q}) = 0$.
\medskip
Note that equalities (of measurable sets or functions) are understood as equalities a.s.
\bigskip

We will suppose that the action $T_{\underline{i}}$ is ergodic, that is, the only sets $A\in \Cal A$ with 
$\mu(A) = \mu(T_{-\underline{i}}A)$ for all $\underline{i} \in \Bbb Z^d$ are of measure 0 or 1.

By the Billingsley-Ibragimov theorem, for $d=1$ (and $\mu$ ergodic) for every martingale 
difference $f\in L^2$ we have a CLT  and the convergence is towards the normal law $\Cal N(0,\sigma^2)$, 
$\sigma^2 = \|f\|_2^2$.
For $d\geq 2$ the central limit theorem is still valid but in general, the limit law is a mixture of normal laws.
More precisely, we have

\proclaim{Theorem 1}  Let $T_{\underline{i}}$ be an ergodic $\Bbb Z^d$ action, let
$(U_{\underline{i}}f)_{\underline{i}}$ (where $U_{\underline{i}}f = f\circ T_{\underline{i}}$) be a field of $L^2$ 
martingale differences for a completely commuting filtration $(\Cal F_{\underline{i}})_{\underline{i}}$. Then for 
$\min(n_1, \dots, n_d) \to \infty$
$$
  \frac1{\sqrt{n_1, \dots, n_d}} \sum_{i_1=1}^{n_1} \ldots \sum_{i_d=1}^{n_d} U_{\underline{i}}f   
  \overset \Cal D \to{\longrightarrow} \nu; \tag{1.1}
$$
$\nu$ is a law with characteristic function $E \exp(-\frac12 \eta^2 t^2)$ where
$$\multline
  \eta^2 =  
  \lim_{n_1, \dots, n_d\to\infty} \frac1{n_1}  \sum_{i=1}^{n_1} \big(\frac1{\sqrt{n_2 \cdots n_d}} \sum_{i_2=1}^{n_2} \ldots \sum_{i_d=1}^{n_d} U_{\underline{i}}f\big)^2 = \\
 \lim_{n_2, \dots, n_d\to\infty} E\Big(   \big(\frac1{\sqrt{n_2 \cdots n_d}} \sum_{i_2=1}^{n_2} \ldots \sum_{i_d=1}^{n_d} U_{\underline{i}}f\big)^2 \,\big| \, \Cal I_1\Big), 
  \endmultline \tag{1.2}
$$
$E\eta^2 = Ef^2$, and $\Cal I_q$ is the sigma algebra of measurable sets invariant for $T_{e_q}$. The convergence 
is in distribution.
\endproclaim

The result was published in \cite{V19}. 
For $d=2$, the result has been carefully proved in \cite{V19} but for $d>2$ more details
are needed.
The present paper thus serves as a kind of an Erratum to \cite{V19}. \newline

If $\eta^2$ is a constant, the limit law is normal. As Wang and Woodroofe showed in their counterexample 
(see \cite{WaWo}),  in general the limit can be different from a normal law.
\medskip

\underbar{\bf Counterexample.} Let $(\Omega_i, \Cal A_i, \mu_i)$, $1\leq i\leq d$, be probability spaces equipped 
with bijective, bimeasurable, and measure preserving mappings $S_i$ of $\Omega_i$. We suppose that all $S_i$
are ergodic and of positive entropy.

On the space $\Omega_1 \times \dots \times \Omega_d$ with the product sigma algebra $\Cal A$ and product
measure $\mu$ we define $T_{\underline{i}}(x_1,\dots,x_d) = (S_1^{i_1}x_1, \dots, S_d^{i_d} x_d)$. The action 
$T_{\underline{i}}$ is ergodic (cf\. \cite{V15}). Like in \cite{GLV} we call it {\it action of product type}.

Because $S_i$ are of positive entropy, there exist ($S_i$-invariant) filtrations $(\Cal F_{j}^i)_j$ and non zero martingale 
differences $f_i\in L^2(\mu_i)$ (defined on $\Omega_i$, with respect to $(\Cal F_{j}^i)_j$; cf\. \cite{V87}). 
On $\Omega$ we define $f((x_1,\dots,x_d)) = f_1(x_1)\cdots f_d(x_d)$. We have $f\in L^2(\mu)$ and $f$
is a martingale difference for the completely commuting filtration $\Cal F_{\underline{i}} = 
\Cal F_{i_1}^1 \otimes \dots \otimes \Cal F_{i_d}^d$. More details are given in \cite{V15} and \cite{GLV}.

Using the Billingsley-Ibragimov CLT we can easily see that
$$
  \frac1{\sqrt{n_1, \dots, n_d}} \sum_{i_1=1}^{n_1} \ldots \sum_{i_d=1}^{n_d} U_{\underline{i}}f =
  \Big(\frac1{\sqrt{n_1}} \sum_{i_1=1}^{n_1} U_1^{i_1} f_1\Big) \cdots
  \Big(\frac1{\sqrt{n_d}} \sum_{i_d=1}^{n_d} U_d^{i_d} f_d\Big)
$$
converge in law to a product of $d$ independent normally distributed random variables hence to a law which is not 
normal.

For convergence towards normal laws, several sufficient conditions are known. In \cite{WaWo} it is shown:

\proclaim{Theorem A (Wang and Woodroofe)}
If the action $T$ is Bernoulli then the limit law in $(1.1)$ is normal.
\endproclaim

The result was improved in

\proclaim{Theorem B (Voln\'y)}
If there is a $1\leq q\leq d$ such that the transformation $T_{e_q}$ is ergodic then the limit law in $(1.1)$ is normal.
\endproclaim

The result was published in \cite{V15}, see also \cite{CDV} where interesting examples are given.
The theorem led to many limit theorems proved by approximation; it turned out that many approximation
results for sequences can be extended to random fields (when using a completely commuting filtration), 
cf\. e.g\. \cite{PZ}, \cite{ZRP}, \cite{G}, \cite{EG}, \cite{CDM}.
 
\subheading{A factor of product type}
For simplicity suppose that $d=2$.
By Theorem B a necessary condition that there is a non normal limit in the CLT for (1.1) is that none of
the transformations $T_{e_1}$, $T_{e_2}$, is ergodic. This means that both sigma algebras $\Cal I_1, \Cal I_2$
are nontrivial. One can notice that the sigma algebras $\Cal I_1, \Cal I_2$ are independent and 
$\Cal I = \Cal I_1\vee \Cal I_2$ is a factor of product type (cf\. \cite{GLV}). This observation started a research presented in \cite{GLV}.

Recall that $\Cal I_q$ denotes the sigma algebra of $T_{e_q}$-invariant sets from $\Cal A$. We define
$$
  \bar \Cal I_q = \bigcap_{1\leq s\leq d, \,s\neq q} \Cal I_s, \quad \Cal I = \bigvee_{1\leq q\leq d} \bar \Cal I_q.
$$

\proclaim{Proposition C (Giraudo, Lesigne, Voln\'y)} The sigma algebras $\bar \Cal I_q$ are factors for the transformations $T_{e_q}$
and are mutually independent. The sigma algebra $\Cal I$ is a factor for the $\Bbb Z^d$ action $T_{\underline{i}}$.

The factor $\Cal I$ is thus isomorphic to the action of product type
$$
  T_{i_1,\dots,i_d}(x_1,\dots,x_d) = (T_{e_1}^{i_1}x_1,\dots, T_{e_d}^{i_d}x_d)
$$
on $\bar \Cal I_1 \otimes \dots \otimes \bar \Cal I_d$ where $T_{e_q}$ acts on $\bar \Cal I_q$.
\endproclaim

We call $\Cal I$ {\it factor of product type}.
\medskip

The Wang-Woodroofe example which we get by the factor $\Cal I$ is a key for understanding convergence to non 
normal distributions.
As we shall see, any $\Bbb Z^d$ action contains a factor of product type and there exists a random field
of martingale differences giving a non normal limit in the CLT if and only if the factor is in some sense non trivial
(see Theorem 2).


\proclaim{Theorem  2} There exists a martingale difference $f\in L^2$ with  non normal limit in $(1.1)$ (i.e\. $\eta^2$ 
is not constant) if and only if all the transformations $T_{e_q}$, $1\leq q\leq d$, are of positive entropy in the 
factor $\Cal I$.
\endproclaim

For $d=2$ the Theorem was proved in \cite{GLV}. 
\medskip

\underbar{Remarks:}

1. We have $\Cal I = (\Cal I \cap \Cal I_q) \vee \bar \Cal I_q$ where $\Cal I \cap \Cal I_q$ and $\bar \Cal I_q$ are 
independent
$T_{e_q}$-invariant sigma algebras hence the entropy of $T_{e_q}$ in $\Cal I$ is the same as in $\bar \Cal I_q$
(it is zero in $\Cal I \cap \Cal I_q$).
By Proposition C, in the Wang-Woodroofe's example we can take $T_{e_q}$ for $S_q$, $\bar \Cal I_q$ for $\Cal A_q$, $\Omega$ for $\Omega_q$. Entropy of  $T_{e_q}$ in $\Cal I$ thus equals entropy of
$S_q$ in $\bar \Cal I_q$. Positive entropy implies existence of a filtration giving a nonzero martingale difference
hence  if all the transformations $T_{e_q}$, $1\leq q\leq d$, are of positive entropy in the factor $\Cal I$ then 
$\Cal I$ is isomorphic to the Wang-Woodroofe's example and we get a martingale difference $f$ for the
$\Bbb Z^d$ action giving convergence to a non normal law.

\comment
2. If the action is Bernoulli then all sigma algebras 
$\Cal I_q$ are trivial hence $\Cal I$ is trivial and by Theorem 2 the limit in the CLT for martingale differences is normal. Theorem 2 thus implies the result in \cite{WaWo}.
\endcomment
2. 
\comment
If one of the transformations $T_{e_q}$ is ergodic (which is a weaker assumption than bernoullicity of the action
and stronger than ergodicity of the action) then in the CLT for martingale differences the limit is normal.
This was proved in [V15], see also [CDV]. Again, the result is a corollary to Theorem 2 above. 
\endcomment
If $T_{e_q}$ is ergodic then $\Cal I_q$ is trivial hence all $\Bar I_r$ with $r\neq q$ are trivial and the transformations 
$T_{e_r}$ are of zero entropy in $\Cal I$. Theorem B is thus a corollary to Theorem 2. Recall that Theorem B implies
Theorem A.

\subheading{A trace of the filtration on $\Cal I$ and a decomposition}

\proclaim{Theorem D (Giraudo, Lesigne, Voln\'y)} Let $f\in L^2$ be a martingale difference adapted to the filtration 
$(\Cal F_{\underline{i}})_{\underline{i}\in \Bbb Z^d}$. The projections $E(f | \Cal I)$ and $f - E(f | \Cal I)$
are martingale differences adapted to the same filtration.
\endproclaim

The theorem is proved in \cite{GLV} but for reader's convenience we show a (different) proof here.

For simplicity, suppose that $d=2$; we denote $T = T_{1,0}$ and $S=T_{0,1}$, we have a completely commuting 
filtration $\Cal F_{i,j}$. The space $L^2$ is decomposed into a direct sum of 
$L^2(\Cal F_{i,\infty}) \ominus L^2(\Cal F_{i-1,\infty})$, $i\in \Bbb Z$, $L^2(\Cal F_{-\infty,\infty})$, and
$L^2  \ominus L^2(\Cal F_{\infty,\infty})$. 
Denote
$$
  P'f = E(f | \Cal I_2) = \lim_{n\to\infty} \frac1n \sum_{i=1}^n f\circ S^i.
$$

\proclaim{Lemma 1} 
$$
  P'(L^2(\Cal F_{i,\infty}) \ominus L^2(\Cal F_{i-1,\infty})) \subset L^2(\Cal F_{i,\infty}) \ominus L^2(\Cal F_{i-1,\infty})
$$
and similarly for $L^2(\Cal F_{-\infty,\infty})$, $L^2  \ominus L^2(\Cal F_{\infty,\infty})$.

Therefore, if $f$ is a martingale difference for the filtration $(\Cal F_{i,\infty})_i$ then $P'f$ is
a martingale difference for the same filtration and also for the filtration $(\Cal I \cap \Cal F_{i,\infty})_i$.
\endproclaim

\demo{Proof}
We notice that for $f\in L^2(\Cal F_{i,\infty}) \ominus L^2(\Cal F_{i-1,\infty})$, $f\circ S^i \in 
L^2(\Cal F_{i,\infty})$ hence $P'f \in L^2(\Cal F_{i,\infty})$. For $g\in L^2(\Cal F_{i-1,\infty})$ we in the same way 
deduce $P'g \in L^2(\Cal F_{i-1,\infty})$ and because $<f\circ S^i, g\circ S^j> =0$ for all $i,j \in \Bbb Z$,
we have $P'f \in L^2(\Cal F_{i,\infty})  \ominus L^2(\Cal F_{i-1,\infty})$.
\enddemo
\qed

Notice that $P'(L^2(\Cal F_{j,\infty}) \ominus L^2(\Cal F_{j-1,\infty}))$ are closed subspaces.

To see this let $f_n \in P'(L^2(\Cal F_{j,\infty}) \ominus L^2(\Cal F_{j-1,\infty}))$ and suppose $f_n \to f$ (in $L^2$).
$P'$ is a surjective mapping of $L^2$ onto $L^2(\Cal I_2)$ hence there is a $h\in L^2$ such that $P'h=f$.
We have $h = \sum_{i\in \Bbb Z \cup \{\infty, -\infty\}} h_i$ where $h_i \in L^2(\Cal F_{i,\infty}) \ominus 
L^2(\Cal F_{i-1,\infty})$ for $i\in \Bbb Z$, $h_{-\infty} \in L^2(\Cal F_{-\infty,\infty})$,
$h_{\infty} \in L^2  \ominus L^2(\Cal F_{\infty,\infty})$. Then $f = P'h = \sum_{i\in \Bbb Z \cup \{\infty, -\infty\}} P'h_i$.
For every $j'\neq j$ it is $\|P'h_{j'}\|_2 \leq \|f_n-f\|_2$ for all $n$ hence $f = P'h = h_j$.
Therefore, 
$P'(L^2(\Cal F_{j,\infty}) \ominus L^2(\Cal F_{j-1,\infty}))$ is a closed subspace of $L^2$ and similarly we show that 
$P'L^2(\Cal F_{-\infty,\infty})$, and $P'(L^2  \ominus L^2(\Cal F_{\infty,\infty}))$ are closed subspaces.
\medskip

Using independence of the sigma algebras $\Cal I_1$, $\Cal I_2$ we can get similar results for the projection
onto $L^2(\Cal I)$. In particular, we can use the fact that if
the family of $e_i$, $i\in I$, is an orthonormal basis of $L^2(\Cal C_1)$, and the family of $f_j$, $j\in J$ is an orthonormal basis of $L^2(\Cal C_2)$, where $\Cal C_1$ is a sub sigma algebra of $\Cal I_1$ and $\Cal C_2$ is a 
sub sigma algebra of $\Cal I_2$, then the family of $e_if_j$, $(i,j) \in I\times J$, is an orthonormal basis of 
$L^2(\Cal C_1\vee \Cal C_2)$.

For general $d\geq 2$ the proof is the same, we use the sigma algebras $\bar \Cal I_q$ (notice that
for $d=2$, $\bar \Cal I_1 = \Cal I_2$ and $\bar \Cal I_2 = \Cal I_1$). 
\bigskip

In \cite{GLV}, all possible limit laws in the CLT for $\Cal I$-measurable martingale differences are given.
\medskip

When starting the research leading to \cite{GLV} I conjectured  that the limit law for $f - E(f | \Cal I)$ is normal 
and the limit law for $f$ is a convolution
of limit laws for $E(f | \Cal I)$ and $f - E(f | \Cal I)$. This is not the case:
in \cite{GLV} a counterexample where the limit law for $f - E(f | \Cal I)$ is not normal is given.

There is, nevertheless, a sufficient condition for normality of the  limit law for $f - E(f | \Cal I)$. 

\proclaim{Proposition 3}
Let $f\in L^2$ be a martingale difference. If for $\ell\to\infty$
$$
  \|E(f\,|\, \Cal F_{-\ell}^{(2)}\vee \Cal I)\|_2 \to 0 \tag{1.3}
$$
then for $\min\{n_1,\dots,n_d\} \to\infty$ 
$$
  \frac1{\sqrt{n_1\cdots n_d}} \sum_{i=1}^m\sum_{i_1=1}^{n_1} \dots \sum_{j_d=1}^{n_d} f\circ T_{i_1,\dots,i_d}
  \to \Cal N(0,\sigma^2) \tag{1.4}
$$
in distribution; $\sigma^2 = \|f\|_2^2$.
\endproclaim

\medskip

\underbar{Remarks.}

1. If $f$ satisfies (1.3) then $E(f | \Cal I) =0$ hence $f\in L^2\ominus L^2(\Cal I)$.

2. For $d=2$ and $f\in L^4$ the Proposition was proved in \cite{GLV}.

\comment
One can easily see that if one of the limit laws for $E(f | \Cal I)$ or for $f - E(f | \Cal I)$ is normal, the limit
law for $f$ is their convolution. 

In the Example in [GLV] there is a random field of $f\circ T_{i,j} = X_iY_jz(i,j)$ where
$(X_i)_i$ and $(Y_j)_j$ are independent sequences of Bernoulli iid, $z(i,j) = (-1)^{i+j}$.
Then  $f\in L^2\ominus L^2(\Cal I)$; for $g \circ T_{i,j} = X_iY_j$ we have $g\in L^2(\Cal I)$.

Let $h = f+g$. The $f,g, h$ are martingale differences.

Then for $\min(n, m) \to\infty$ the limit law for $(1/\sqrt{nm})S_{n,m}(f)$ as well as for $(1/\sqrt{nm})S_{n,m}(g)$   
is the product of two independent r.v\. with $\Cal N(0,1)$ distribution.

In the sum $(1/\sqrt{nm})S_{n,m}(f+g)$ we have terms $V_{i,j} = X_iY_j (z(i,j)+1)$ hence 
the sum of $V_{i,j}$ is the sums of $2X_iY_j$ where $i, j$ are both even plus the sum where $i, j$ are both odd.
Notice that the "even" sum is independent of the "odd" sum. $(1/\sqrt{nm})S_{n,m}(f+g)$ thus converge (in 
distribution) towards a convolution of two normal laws, hence a normal law. 
On the other hand, the convolution of the limits of $(1/\sqrt{nm})S_{n,m}(f)$ and $(1/\sqrt{nm})S_{n,m}(g)$
is not a normal law.

The Example can be modified in a way that as soon as there exists a martingale difference $f\in L^2$ leading
to a non normal limit in the CLT, we can have martingale differences $f'\in L^2(\Cal I)$ and $f''\in 
L^2\ominus L^2(\Cal I)$ such that both $f'$ and $f''$ give non normal limit laws $\nu', \nu''$ in the CLT and the limit
law for $f'+f"$ is not the convolution $\nu'*\nu''$.

Determination of all possible limit laws in the CLT for any martingale difference $f\in L^2$ and for $f - E(f | \Cal I)$
remains an open problem. 
\endcomment

\subheading{An open problem}

The question what is the limit law (what function $\eta^2$ we can have) for a martingale difference $f$ and
for a martingale difference $f\in L^2\ominus L^2(\Cal I)$ in particular, remains open.

\heading
2. Proof of Theorem 1
\endheading

The proofs use McLeish's martingale CLT \cite{McL}.

\proclaim{Theorem E (McLeish)} Let $X_{n,k}$, $1\leq k\leq k_n$, be an array of martingale differences such that
\roster
\item"(i)" $\max_{1\leq k \leq k_n} |X_{n,k}|$ are uniformly bounded in $L^2$ norm,
\item"(ii)" $\max_{1\leq k \leq k_n} |X_{n,k}| \to 0$ in probability,
\item"(iii)" $\sum_{k=1}^{k_n} X_{n,k}^2$ converge in probability to a constant $\sigma^2=0$.
\endroster
Then $\sum_{k=1}^{k_n} X_{n,k}$ converge in distribution to the normal law $\Cal N(0,\sigma^2)$.
\endproclaim

In our case we take martingale differences $f\circ T^i$, $X_{n,k} = f\circ T^k/\sqrt n$ ($1\leq k\leq n$).
It is well known that for martingale differences $f\circ T^i \in L^2$, (i) and (ii) are satisfied.
In (iii) we by Ergodic Theorem get $(1/n) \sum_{i=1}^n f^2\circ T^i \to \eta^2 = E(f^2 | \Cal I)$ where
$\Cal I$ is the sigma algebra of $T$-invariant measurable sets; the convergence is in $L^1$ and almost surely.

Without loss of generality we can suppose that there exists a regular family of conditional probabilities with respect
to $\Cal I$ (cf\. \cite{V89}). The conditional probabilities are ergodic and the martingale property is preserved. 
For a.e\. 
conditional probability we thus can use McLeish's theorem and eventually we get a convergence towards a mixture
of normal laws (cf\. \cite{HaHe} and \cite{V89}).

\proclaim{Proposition F} The random variables $(1/\sqrt n)  \sum_{i=1}^n f \circ T^i$ converge in distribution
to a law with characteristic function $\varphi(t) = E \exp(-\frac12 \eta^2t^2)$; $\eta^2 = E(f^2 | \Cal I)$.
\endproclaim

Recall that if $T$ is ergodic then $\eta^2 = E(f^2 | \Cal I)$ is a constant.
\medskip

\subheading{Proof of Theorem 1}

For $d=2$, a proof (detailed and correct) can be found in \cite{V19}. In the same paper ideas how to extend it to 
$d\geq 3$ are given but more arguments are needed.
In particular, in the proof of Lemma 8 it is shown that $E\eta_1^4 \geq E \eta_2^4$ and by symmetry 
$E\eta_2^4 \geq E \eta_1^4$ hence we get an  equality.
Within the framework of the proof for $d=2$ in \cite{V19} the symmetry is evident; this is not the case for $d>2$. 
Here we present a proof for $d\geq 3$ using results from \cite{V19}.

From now on we suppose $d\geq 3$ and that the theorem is true for $d-1$.
\medskip

We denote $\underline{v} = (n_2,\dots, n_d)$ and define
$$
  F_{i,\underline{v}} = F_{ i,n_2,\dots,n_d} = 
  \frac1{\sqrt{n_2\cdots n_d}} \sum_{i_2=1}^{n_2} \dots \sum_{i_d=1}^{n_d} U_{i,i_2\dots,i_d}f.
$$

\proclaim{Proposition 4} If for $n_1\to\infty$
there exist $n_2(n_1), \dots n_d(n_1)$ such that for \newline
$n_2\geq n_2(n_1), \dots n_d\geq n_d(n_1)$ uniformly
$$
  \Big\| \frac1{n_1}  \sum_{i=1}^{n_1} F_{ i,n_2,\dots,n_d}^2 - 
   E\big( F_{1,n_2,\dots,n_d}^2 \,  \big|\, \Cal I_1\big) \Big\|_1 < \epsilon(n_1) \to 0  \tag{2.1}
$$
and 
$$
  E\big( F_{1,n_2,\dots,n_d}^2 \,  \big|\, \Cal I_1\big) \overset \Cal D \to{\longrightarrow} \eta^2  \quad
  \text{as}\quad \min(n_2,\dots,n_d) \to \infty \tag{2.2}
$$
then $(1.1)$ holds true.
\endproclaim

\demo{Proof} In the same way as in the proof of Theorem 1 in \cite{V19}, from (1.1) can be deduced that for any 
sequence
of $\underline{v}_n \to \underline{\infty}$ and any sequence of $N_n$ with $N_n/n\to \infty$ there exist
$\Cal I_1$-measurable random variables $\eta^2(n) \overset \Cal D \to{\longrightarrow} \eta^2$ such that 
$$
  \Big\| \frac1{N}  \sum_{i=1}^{N} F_{i,\underline{v}}^2 - \eta(n)^2 \Big\|_1 \to 0.
$$
From this, in the same way as in the proof of Theorem 1 in \cite{V19} we can see that assumptions of 
Proposition 3 in \cite{V19} are satisfied.
\enddemo
\qed

\demo{Proof of Theorem 1} 
In \cite{V19} the theorem has been proved for $d=2$. Recall that here we suppose $d\geq 3$ and the theorem
proved for $d-1$.
We thus suppose $d\geq 3$.

We can suppose that $f$ is bounded: Let $P_{\underline{0}}$ be the orthogonal projection on
$\big( L^2(\Cal F_{\underline{0}})\ominus L^2(\Cal F_{-e_1}) \big) \cap \dots \cap 
\big( L^2(\Cal F_{\underline{0}})\ominus L^2(\Cal F_{-e_d}) \big)$ (recall that $e_q$ is the vector $(j_1,\dots,j_d)
\in \Bbb Z^d$ where $j_q=1$ and $j_i=0$ for all $i\neq q$). For a $K>0$, $ P_{\underline{0}}(f 1_{|f|\leq K})$
and $f - P_{\underline{0}}(f 1_{|f|\leq K})$ are martingale differences. For $K\to \infty$ we have
$\|f - P_{\underline{0}}(f 1_{|f|\leq K})\|_2 \to 0$ hence we can replace $f$ by a bounded function 
$P_{\underline{0}}(f 1_{|f|\leq K})$ (cf\. \cite{V19}).

\comment
By validity of the theorem for $d-1$, for $\min(n_2\dots,n_d) \to \infty$ we get
$$
  (F_{i,n_2\dots,n_d})_i \overset \Cal D \to{\longrightarrow} (V_{1,i})_i
$$
where $(V_{1,i})_i$ is a (stationary) sequence of square integrable martingale differences and $\|V_{1,i}\|_2^2 = 
\|f\|_2^2$, $V_{1,i}\in L^p$ for all $1\leq p<\infty$. This can be proved in the same way as Lemma 5 in \cite{V19}.
More precisely, for any finite set $J\subset \Bbb Z$ and real numbers $a_i, i\in J$, $\sum_{i\in J} a_i U_{i,\underline{j}}$
with $\underline{j} \in \Bbb Z^{d-1}$ are martingale differences hence $\sum_{i\in J} a_i F_{i,\underline{j}}$ converge 
in distribution.
\endcomment

We denote $\underline{v} = (n_2,\dots, n_d)$ and define
$$
  F_{i,\underline{v}} = F_{ i,n_2,\dots,n_d} = 
  \frac1{\sqrt{n_2\cdots n_d}} \sum_{i_2=1}^{n_2} \dots \sum_{i_d=1}^{n_d} U_{i,i_2\dots,i_d}f.
$$

\proclaim{Lemma 2}
For $\min(n_2\dots,n_d) \to \infty$
$$
  (F_{i,n_2\dots,n_d})_i \overset \Cal D \to{\longrightarrow} (V_{1,i})_i
$$
where $(V_{1,i})_i$ is a (stationary) sequence of square integrable martingale differences and $\|V_{1,i}\|_2^2 = 
\|f\|_2^2$, $V_{1,i}\in L^p$ for all $1\leq p<\infty$. 
\endproclaim

\demo{Proof of the Lemma}

To prove this, take a finite $J\subset \Bbb Z$.
For any vector 
$(a_i: i\in J)\in \Bbb R^J$, $\sum_{i\in J} a_i U_{i,0,\dots,0}f$ is a martingale difference for $\Bbb Z^{d-1}$ action of
$T_{0,i_2,\dots,i_d}$ hence by assumption of validity of Theorem 1 for $d-1$
$$
  \frac1{\sqrt{n_2\cdots n_d}} \sum_{i_2=1}^{n_2} \dots \sum_{i_d=1}^{n_d} \sum_{i\in J} a_i U_{i,i_2\dots,i_d}f
$$
converge in distribution to a law with characteristic function $\exp(-\frac12 \eta^2t^2)$ where
$$
  \eta^2 = \lim_{\min(n_3,\dots,n_d)\to\infty} E\big( (\sum_{i\in J} a_i F_{i,1,n_3,\dots,n_d})^2 \,\big|\, \Cal I_2\big).
$$
The probability of the event that $(F_{i, \underline{n}} : i\in J)$ is not
in $[-K, K]^J$ thus goes with $K\to \infty$ to zero uniformly for all $\underline{n}$ hence the set of (distributions of) 
$(F_{i, \underline{n}} : i\in J)$ is tight. By the Cramer-Wold device there is a unique limit point, therefore the net
$(F_{i, \underline{n}} : i\in J)$ converge.

The family of finite dimensional distributions of $(V_{1,i})_{i\in J}$ is projective hence by the Kolmogorov's 
theorem we get a distribution $\nu_1$ of $(V_{1,i})_{i\in \Bbb Z}$.
\bigskip

From stationarity of $(F_{i,\underline{v}})_i$ we get stationarity of  $(V_{1,i})_i$.
Because $\|F_{i,\underline{v}}\|_2 = \|f\|_2$, using \cite{B68, Theorem 5.3} (cf\. (2.15) here) we can show
$\|V_{1,i}\|_2 = \|f\|_2$. 

The characteristic function of $V_{1,1}$ is $E \exp(-\frac12 \eta^2t^2)$ where 
$$
  \eta^2 = \lim_{\min(n_3,\dots,n_d)\to\infty} E( F_{1,1,n_3,\dots,n_d}^2 | \Cal I_2).
$$
By validity of Theorem 1 for $d-1$ we deduce finite $p$-th moments of $\eta^2$ hence of $V_{1,1}$ for all 
$1\leq p<\infty$.

In the same way as in \cite{V19, Lemma 5} we show that $(V_{1,i})_i$ is a sequence of martingale differences.
It is sufficient to prove that for any bounded measurable function $g$ on $\Bbb R^n$ and $n\geq 1$,
$\int V_{1,0} g(V_{1,-n},\dots,V_{1,-1}) \,d\nu_1 =0$; we have 
$0 = E\Big(F_{0,\underline{v}}g\big(F_{-n,\underline{v}},\dots,F_{-1,\underline{v}}\big)\Big) \to 
\int V_{1,0} g(V_{1,-n},\dots,V_{1,-1}) \,d\nu_1$ as $\min(v_1,\dots,v_d) \to \infty$.

\enddemo
\qed
\bigskip

\underbar{Remark.} The Lemma is true for $d=2$ as well. It was proved in \cite{V19, Lemma 5}. The present proof
remains valid; Theorem 1 then becomes Proposition F.
\medskip

By  Ergodic Theorem there exists a function $\eta_1^2$ such that
$$
  \frac1n  \sum_{j=1}^n V_{1,j}^2 \to  \eta_1^2 \tag{2.3}
$$
(in all $L^p$,  $1\leq p<\infty$) and by Proposition F $(1/\sqrt n) \sum_{j=1}^n V_{1,j}$ converge in distribution to 
a law with characteristic function $\exp(-\frac12 \eta_1^2 t^2)$. 

\bigskip

For $(n_2, \dots n_d)\in \Bbb N^{d-1}$ fixed we (using Ergodic Theorem) define
$$
  \eta_{1,n_2, \dots n_d}^2 = \lim_{n_1\to\infty} \frac1{n_1}  \sum_{i=1}^{n_1} F_{i,n_2, \dots, n_d}^2 =
  E(F_{1,n_2, \dots, n_d}^2 | \Cal I_1); \tag{2.4}
$$
$F_{ i,n_2, \dots, n_d}$ are martingale differences and by Proposition F the random variables \newline
$(1/\sqrt{ n_1}) 
\sum_{i=1}^{n_1} F_{i, n_2, \dots, n_d}$ converge in distribution to a law with characteristic function  
$\exp(-\frac12 \eta_{1, n_2, \dots, n_d}^2 t^2)$. 

\bigskip

In the same way as before we can define 
$$
  F_{n_1, j,n_3,\dots,n_d} = \frac1{\sqrt{n_1 n_3\cdots n_d}} \sum_{i_1=1}^{n_1} \sum_{i_3=1}^{n_3}
  \dots \sum_{i_d=1}^{n_d} U_{i_1,j,i_3,\dots,i_d}f;
$$
 for $\min(n_1, n_3,\dots,n_d) \to \infty$ we get
$$
  (F_{n_1, j,n_3,\dots,n_d})_j \overset \Cal D \to{\longrightarrow} (V_{2,j})_j
$$
(where $(V_{2,j})_j$ is a strictly stationary sequence of martingale differences) and similarly as in (2.3)
$$
  \frac1n  \sum_{j=1}^n V_{2,j}^2 \to  \eta_2^2
$$
(in $L^1$ and a.s.). By Proposition F  $(1/\sqrt n) \sum_{j=1}^n V_{2,j}$ converge in distribution to a law with 
characteristic function  $E\exp(-\frac12 \eta_2^2 t^2)$.

We can notice that
$$
  \frac1{\sqrt{n_1}} \sum_{i=1}^{n_1}  F_{i,n_2,\dots,n_d} =  
  \frac1{\sqrt{n_2}} \sum_{j=1}^{n_2}  F_{n_1, j,n_3,\dots,n_d}.
$$
Recall that for $n_2, \dots, n_d$ fixed and $n_1\to \infty$ the first sums converge in distribution to a law with 
characteristic function $E\exp(-\frac12 \eta_{1,n_2,\dots,n_d}^2 t^2)$, and by validity of the theorem for $d-1$
there exist $n_1(n_2), n_3(n_2), \dots, n_d(n_2)$ and a function 
$\eta_2^2$ such that for $n_2\to\infty$ and $n_1\geq n_1(n_2), n_3\geq n_3(n_2), \dots, n_d\geq n_d(n_2)$
uniformly, the second sums converge in distribution to a law with characteristic function 
$E\exp(-\frac12 \eta_2^2 t^2)$.

From this we deduce that for $n_2\to\infty$ and for $n_3\geq n_3(n_2),\dots, n_d\geq n_d(n_2)$ uniformly,
$$
  E \exp(-\frac12 \eta_{1,n_2,\dots,n_d}^2 t^2) \to E  \exp(-\frac12 \eta_2^2 t^2)
$$
for all $t$ hence by properties of Laplace transformation 
$$
  \eta_{1,n_2,n_3,\dots, n_d}^2 \overset \Cal D \to{\longrightarrow} \eta_2^2 \tag{2.5}
$$
(see \cite{F}, Chapter XIII).

In the same way as in Lemma 2 and the Remark after it
(cf\. also proof of Lemma 5 in \cite{V19})
we can prove that for $(j_2,\dots, j_d) \in \Bbb Z^{d-1}$ fixed and $n_1\to\infty$, 
$$
  \tilde F_{n_1, j_2, \dots, j_d} = \frac1{\sqrt{n_1}} \sum_{i_1=1}^{n_1} U_{i_1,j_2\dots,j_d}f 
  \overset \Cal D \to{\longrightarrow} W_{j_2,\dots, j_d}
  \tag{2.6}
$$
where $(W_{j_2,\dots, j_d})_{(j_2,\dots, j_d)\in \Bbb Z^{d-1}}$ is a field of martingale differences.
By validity of the theorem for $d-1$ there exists a function $\eta_{2,\dots,d}^2$ such that for
$\min(n_2, \dots, n_d) \to \infty$ the sums 

$(1/\sqrt{n_2 \cdots n_d}) \sum_{j_2=1}^{n_2} \dots \sum_{j_d=1}^{n_d} W_{j_2,\dots, j_d}$ \newline
converge (in distribution) to a law with  characteristic function $E \exp(-\frac12 \eta_{2,\dots,d}^2 t^2)$ and 
$$
  \frac1{n_2} \sum_{j_2=1}^{n_2} \Big( \frac1{\sqrt{n_3\cdots n_d}} \sum_{j_3=1}^{n_3} \dots \sum_{j_d=1}^{n_d}
  W_{j_2,\dots,j_d} \Big)^2 \overset \Cal D \to{\longrightarrow}  \eta_{2,\dots,d}^2.
$$
There thus exist $n_1( n_2), \dots, n_1( n_d)$ such that
if $\min(n_2, \dots, n_d) \to \infty$ and $n_1 \geq \max( n_1( n_2), \dots, n_1(n_d))$ then
$$
  \frac1{n_1 \cdots n_d} \sum_{i_2=1}^{n_2} \Big(  \sum_{i_1=1}^{n_1}  \sum_{i_3=1}^{n_3} \dots \sum_{i_d=1}^{n_d} 
  U_{i_1,\dots,i_d}f  \Big)^2 \overset \Cal D \to{\longrightarrow}  \eta_{2,\dots,d}^2.
$$

If $n_2,  \dots, n_d$ are big enough and $n_1 \geq \max(n_1( n_2), \dots, n_1(n_d))$ then  the distribution 
of
$$
  \frac1{\sqrt{n_1, \dots, n_d}} \sum_{i_1=1}^{n_1} \dots \sum_{i_d=1}^{n_d} U_{i_1,\dots,i_d}f 
$$
is close to a law  with  characteristic function $E \exp(-\frac12 \eta_{2, \dots, d}^2 t^2)$ and by the text after (2.4) 
it is also close 
to a law  with  characteristic function $E \exp(-\frac12 \eta_{1,n_2,\dots,n_d}^2 t^2)$. \newline
In the same way as in proving (2.5) (see also the proof of Lemma 7 in \cite{V19}) we deduce that 
$$
  \eta_{1,n_2,\dots,n_d}^2   \overset \Cal D \to{\longrightarrow}  \eta_{2,\dots,d}^2 \tag{2.7}
$$ 
(for $\min(n_2, \dots, n_d) \to\infty$). \newline
By (2.5) for $n_2\to\infty$ and for $n_3\geq n_3(n_2),\dots, n_d\geq n_d(n_2)$ uniformly, 
$ \eta_{1,n_2,\dots,n_d}^2  \overset \Cal D \to{\longrightarrow} \eta_2^2$.
Therefore, 
$$
   \eta_{2,\dots,d}^2 \overset \Cal D \to{=} \eta_2^2. \tag{2.8}
$$
By symmetry we get $\eta_{2,\dots,d}^2 \overset \Cal D \to{=} \eta_3^2$ and in this way we can deduce
$$
   \eta_1^2 \overset \Cal D \to{=} \eta_2^2 \overset \Cal D \to{=} \dots \overset \Cal D \to{=} \eta_d^2. \tag{2.9}
$$
\bigskip

We can (cf\. \cite{V19}, (5) there) choose a sequence of $v_{2,n}, \dots, v_{d,n}$ such that for $n\to\infty$
$$
  Y_n^2 = \frac1{n}  \sum_{i=1}^{n} F_{ i, v_{2,n}, \dots, v_{d,n}}^2 \overset \Cal D \to{\longrightarrow} \eta_1^2 
  \tag{2.10}
$$
(cf\. (2.3); this convergence holds true as soon as $\min(v_{2,n}, \dots, v_{d,n})$ is big enough, say all terms
bigger than a $v(n) \to \infty$).

\proclaim{Lemma 3}
$$
  Y_n^2 - E(Y_n^2 | \Cal I_1) \longrightarrow 0\quad\text{in}\quad L^1.
$$
\endproclaim

\demo{Proof} The proof is very similar to that
of Lemma 8 in \cite{V19} with a few differences caused by the fact that here, $d>2$. 
\medskip
Let us present it in detail. We denote $\underline{v}_n = (v_{2,n}, \dots, v_{d,n})$.

\bigskip

For $k\geq 0$ define
$$
  \varphi_k(x) = \cases  x\,\,\,&\text{if}\,\,\, |x| \leq k, \\
                                      k.\operatorname{sign}(x) &\text{if}\,\,\, |x|>k.
                          \endcases 
$$
By Jensen's inequality taken conditionally we have, for every $k, n\geq 1$,
$$
  E\Big[\varphi_k(Y_n^2) \Big]^2 \geq E\Big[ E\big( \varphi_k(Y_n^2)\,|\,\Cal I_1\big) \Big]^2.  \tag{2.11}
$$
By (2.10),
$
   E\Big[\varphi_k(Y_n^2) \Big]^2 \underset n\rightarrow \infty \to{\longrightarrow} E\Big[\varphi_k(\eta_1^2) \Big]^2
$
for every $k$. We thus have
$$
  \lim_{k\to\infty} \lim_{n\to \infty} E\Big[\varphi_k(Y_n^2) \Big]^2  =  E \eta_1^4.   \tag{2.12}
$$
By concavity of $\varphi_k$ (on $[0,\infty)$) we have  
$$
  (1/n) \sum_{i=1}^n \varphi_k(F_{i,\underline{v}_n}^2) \leq \varphi_k\big( (1/n)  \sum_{i=1}^n 
  F_{i,\underline{v}_n}^2 \big) =   \varphi_k(Y_n^2);
$$ 
we deduce that for every $k$,
$$
   E\big(\varphi_k(F_{1,\underline{v}_n}^2)  \,|\, \Cal I_1\big) = E\Big(\frac1n \sum_{i=1}^n 
  \varphi_k(F_{i,\underline{v}_n}^2)  \,|\, \Cal I_1\Big) \leq 
   E\big(\varphi_k(Y_n^2) \,|\, \Cal I_1\big).  \tag{2.13}
$$

Recall that for $n\to\infty$
$$
 F_{1,\underline{v}_n} \overset \Cal D \to{\longrightarrow} V_1,   \,\,\,\,\,
  E  F_{1,\underline{v}_n}^2 =E f^2 = E V_1^2,\,\,\,\,\, 
$$
and by (2.4), (2.7), (2.8),
$$
  E(F_{1,\underline{v}_n}^2\, |\, \Cal I_1)  \overset \Cal D \to{\longrightarrow} \eta_2^2.
$$
Because $F_{1,\underline{v}_n} \overset \Cal D \to{\longrightarrow} V_{1,1}$ where $\|V_{1,1}\|_2^2 = \|f\|_2^2$
and $E F_{1,\underline{v}_n}^2 = \|f\|_2^2$ for all $n$ we get that the $F_{1,\underline{v}_n}^2$ are uniformly 
integrable. 
Hence, for every $\epsilon>0$ 
there exist 
$k\geq 1$ and $n(k)$ such that for all $n\geq n(k)$,
$$
  E \big| (F_{1,\underline{v}_n}^2 - \varphi_k(F_{1,\underline{v}_n}^2)\big| < \epsilon.
$$
For a given $\delta>0$ we therefore can choose $n$ big enough so that
$$
  \mu\big\{ |  E(F_{1,\underline{v}_n}^2\, |\, \Cal I_1)  - E( \varphi_k(F_{1,\underline{v}_n}^2)  \,|\, \Cal I_1)| >\delta \big\} 
  <\delta. \tag{2.14}
$$
The next statement follows from \cite{B68, Theorem 5.3}.
Let $X, X_n$, $n\geq 1$, be random variables, $p\geq 1$.
$$
  \text{If}\,\,\,\,  X_n \overset \Cal D \to{\longrightarrow} X \,\,\,\, \text{then}\,\,\,\,  E |X|^p \leq  \liminf_{n\to\infty} E |X_n|^p. \tag{2.15}
$$
From $\big[ E(F_{1,\underline{v}_n}^2  \,|\, \Cal I_1)\big]^2  \overset \Cal D \to{\longrightarrow} \eta_2^4$ and (2.15) it follows
$$
   \liminf_{n\to\infty} E[E(F_{1,\underline{v}_n}^2\, |\, \Cal I_1)]^2 \geq E \eta_2^4. \tag{2.16}
$$
For all $1\leq p<\infty$ we have $\eta_1^2\in L^p$, hence by (2.9) $ E \eta_2^4 <\infty$.
By (2.14), (2.16), and uniform integrability of $ \eta_2^4$, for every $\epsilon>0$ there are $k(\epsilon)$ and $n(\epsilon, k)$ such that for all 
$k\geq k(\epsilon), n\geq n(\epsilon, k)$, 
$$
    E\Big[E\big(\varphi_k(F_{1,\underline{v}_n}^2)  \,|\, \Cal I_1\big) \Big]^2 \geq  E \eta_2^4 -\epsilon.
$$
Therefore,
$$
  \lim_{k\to\infty} \liminf_{n\to\infty}  E\Big[E\big(\varphi_k(F_{1,\underline{v}_n}^2)  \,|\, \Cal I_1\big) \Big]^2 \geq  E \eta_2^4. \tag{2.17}
$$
From (2.12), (2.11), (2.13), (2.17) (recall (2.9)) we deduce that for every $\epsilon>0$ there are $k(\epsilon)$ and $n(k,\epsilon)$ such that for $k\geq k(\epsilon)$  and 
$n\geq n(k,\epsilon)$ 
$$\multline
  E \eta_1^4 +\epsilon \geq  E\big[\varphi_k(Y_n^2) \big]^2 \geq  E\Big[ E\big( \varphi_k(Y_n^2) \,|\, \Cal I_1 \big) 
  \Big]^2 \geq \\
  E\Big[ E( \varphi_k(F_{1,\underline{v}_n}^2 \,|\, \Cal I_1 \big) \Big]^2 \geq E \eta_1^4 -\epsilon. 
  \endmultline \tag{2.18}
$$

Using (2.18) we deduce
$$  
   \lim_{k\to\infty} \limsup_{n\to \infty} \Big\{ E\big[\varphi_k(Y_n^2) \big]^2 - E\big[ E\big( \varphi_k(Y_n^2)\,|\,\Cal I_1\big) \big]^2 \Big\} = 0.
$$
This implies
$$
   \lim_{k\to\infty} \limsup_{n\to \infty} E\big[\varphi_k(Y_n^2) -  E\big( \varphi_k(Y_n^2) \,|\, \Cal I_1 \big) \big]^2 = 0,
$$
hence
$$
   \lim_{k\to\infty} \limsup_{n\to \infty} \big\|\varphi_k(Y_n^2) -  
  E\big( \varphi_k(Y_n^2) \,|\, \Cal I_1 \big) \big\|_1 = 0. \tag{2.19}
$$
Because $E Y_n^2 = E f^2 = E \eta_1^2$ and $Y_n^2 \overset \Cal D \to{\longrightarrow} \eta_1^2$, $Y_n^2$ are 
uniformly integrable; we thus have
$$
   \lim_{k\to\infty} \limsup_{n\to \infty}  \| Y_n^2 - \varphi_k(Y_n^2) \|_1 = 0.
$$
From this and (2.19) we deduce
$$
   \lim_{k\to\infty} \limsup_{n\to \infty}  \big\| Y_n^2 -  E\big( \varphi_k(Y_n^2) \,|\, \Cal I_1 \big) \big\|_1 = 0, 
$$
hence
$$
   \lim_{k\to\infty} \limsup_{n\to \infty}  \big\|  E\big( Y_n^2 \,|\, \Cal I_1 \big) -  
  E\big( \varphi_k(Y_n^2) \,|\, \Cal I_1 \big)\big\|_1 = 0,
$$ 
hence
$$
   \lim_{n\to \infty} \big\|  Y_n^2 -  E\big( Y_n^2 \,|\, \Cal I_1 \big) )\big\|_1 = 0;
$$
this finishes the proof of the Lemma.

\enddemo
\qed
\bigskip

Lemma 3 implies (2.1). 

More precisely, by the Lemma (cf\. the proof of Theorem 1 in \cite{V19}) we can show that for 
$n_1\to\infty$ there exist $n_2(n_1), \dots , n_d(n_1)$ and  $\epsilon(n_1) \to 0$ such that for \newline
$n_2\geq n_2(n_1), \dots , n_d\geq n_d(n_1)$
$$
  \Big\| \frac1{n_1}  \sum_{i=1}^{n_1} F_{ i,n_2,\dots,n_d}^2 - 
   E\big( F_{i,n_2,\dots,n_d}^2 \,  \big|\, \Cal I_1\big) ^2 \Big\|_1 <   \epsilon(n_1).   \tag{2.1}
$$
By (2.4),  (2.7), (2.8), (2.9) we have
$$
  E\Big( \big(\frac1{\sqrt{n_2\cdots n_d}} \sum_{i_2=1}^{n_2} \dots \sum_{i_d=1}^{n_d} U_{1,i_2,\dots i_d}f\big)^2 \,
  \big|\, \Cal I_1\Big)  \overset \Cal D \to{\longrightarrow} \eta_1^2
$$
where $\eta_1^2 = \eta^2$.
The convergence in the CLT now follows by Proposition 4, similarly as in the proof of Theorem 1 in \cite{V19}.
\enddemo
\qed

\heading
3. Proofs of Proposition 3 and Theorem 2
\endheading

We first prove a special case of Proposition 3 for $d=2$.

\proclaim{Proposition 3a}
Let $f\in L^2$ be a martingale difference. 
If for $\ell\to\infty$
$$
  \|E(f\,|\, \Cal F_{\infty, -\ell}\vee \Cal I_1)\|_2 \to 0 \tag{3.1}
$$
then for
$$
  F_{i,n} = \frac1{\sqrt n} \sum_{j=0}^{n-1}  f \circ T_{i,j}, 
$$
$$
  \lim_{n\to\infty} \lim_{m\to\infty} \frac1m \sum_{i=0}^{m-1} F_{i,n}^2 = \lim_{n\to\infty} E(F_{1,n}^2\,|\, \Cal I_1)
  = Ef^2  \tag{3.2}
$$
and the convergence is in $L^1$.
\endproclaim


\demo{Proof of Proposition 3a}

First, we truncate $f$ in a more subtle way than in the proof of Theorem 1.
\bigskip
\proclaim{Lemma 3}
Let $\epsilon_k\searrow 0$, $\eta_{1,k}\searrow 0$. For every $\epsilon_k, \eta_{1,k}$ there exist $C'_k$ and $n_k$ such that
for all $n\geq n_k$ and $f'_k =  f 1_{|f|\leq C'_k} - E(f 1_{|f|\leq C'_k} | \Cal F_{\infty,-1})$, $f''_k = f - f'_k$,
$F'_{k;i,n} = (1/n) \sum_{j=0}^{n-1}  f'_k  \circ T_{i,j}$, $F''_{k;i,n} = (1/n) \sum_{j=0}^{n-1}  f''_k  \circ T_{i,j}$,
we have
$$\gathered
  \|E(F_{i,n} | \Cal I_1) - E(F'_{k;i,n} | \Cal I_1)\|_2^2 <  \eta_{1,k}, \\
  \Big\| \frac1n \sum_{j=0}^{n-1}  E({f'_k}^2| \Cal I_1) \circ T_{0,j} - E{f'_k}^2 \Big\|_1 < \epsilon_k.
  \endgathered
   \tag{3.3}
$$
\endproclaim

\demo{Proof} Let us omit the index $k$.
For a $C'>0$ we have decomposed $f$ into a sum of two martingale differences $f=f'+f''$
(we remark that both $f'$ and $f''$ are martingale differences for $T_{0,1}$).
Therefore
$$
  \|F'_{i,n}\|_2^2 = \|f'\|_2^2, \quad \|F''_{i,n}\|_2^2 = \|f''\|_2^2.
$$
Choosing $C'$ big enough, for a given $\eta_1>0$, we can have 
$$
  \|f''\|_2 < \eta_1, \quad  \|f\|_2 \geq \|f'\|_2 >  \|f\|_2 - \eta_1. \tag{3.4}
$$
By Lemma 1 (see also \cite{GLV})
$E(f \circ T_{0,j} \Big| \Cal I_1)$, $E(f' \circ T_{0,j} \Big| \Cal I_1)$, and $E(f'' \circ T_{0,j} \Big| \Cal I_1)$
are martingale differences hence
$$
  \|E(F_{0,n} | \Cal I_1) - E(F'_{0,n} | \Cal I_1)\|_2^2 = \|E(F''_{0,n} | \Cal I_1)\|_2^2 \leq \|f''\|_2 < \eta_1.  \tag{3.5}
$$
This proves the first inequality in (3.3). 

To prove the second inequality we note that by the Ergodic Theorem and independence of $\Cal I_1$ and $\Cal I_2$ 
$$
   \frac1n \sum_{j=0}^{n-1}  E({f'_k}^2| \Cal I_1) \circ T_{0,j} \to E\big( E({f'_k}^2| \Cal I_1) \,\big| \,\Cal I_2\big)
  = E{f'_k}^2.
$$

\enddemo
\qed

We  suppose that the sequence $(C'_k)_k$ is nondecreasing. \newline
Remark that making $n$ bigger the inequalities of (3.5) remain valid.

Now, for $n_k\leq n< n_{k+1}$ we define 
$$
  C_n = C_k', \quad f_n = f'_k =  f 1_{|f|\leq C_n} - E(f 1_{|f|\leq C_n} | \Cal F_{\infty,-1}), \quad 
  \bar F_{i,n} = \frac1{\sqrt n} \sum_{j=0}^{n-1}  f_n \circ T_{i,j}.
$$

For $n\to \infty$ and $\ell_n\to \infty$ we  have 

\proclaim{Lemma 4}
$$\gathered
  \|E(F_{i,n} | \Cal I_1) - E(\bar F_{i,n} | \Cal I_1)\|_2 \to 0, \\
  \Big\| \frac1n \sum_{j=0}^{n-1} E( f_n^2 | \Cal I_1) \circ T_{0,j} - Ef_n^2 \Big\|_1 \to 0, \\
  |Ef^2 - Ef_n^2| \to 0, \\
  \| E( f_n | \Cal F_{\infty,-\ell_n} \vee \Cal I_1) \|_2 \to 0.
  \endgathered \tag{3.6}
$$
\endproclaim

The Lemma follows from (3.3), (3.4) and (3.1) ((3.1) together with the third statement of the Lemma implies the fourth 
statement).

Notice that $(C_n)_n$ is a nondecreasing sequence (because $(C'_k)_k$ is nondecreasing) and 
$\lim_{n\to\infty} C_n = \lim_{k\to\infty} C'_k$. If $f$ is unbounded then $\lim_{n\to\infty} C_n = \infty$.  

If we replace the sequence of $C_n$ by $\bar C_n<C_n$ with the same limit,
(3.6) remains true. 

\medskip

By (3.6) $\|E(F_{i,n} | \Cal I_1) - E(\bar F_{i,n} | \Cal I_1)\|_2 \to 0$ and by (3.4) $| E f^2 - E{f_n}^2|\to 0$.
Because by  triangular inequality
$$\multline
  | E(F_{1,n}^2\,|\,\Cal I_1) - Ef^2| \leq \\
  | E(F_{1,n}^2\,|\,\Cal I_1) - E(\bar F_{1,n}^2\,|\,\Cal I_1) | +
  | E(\bar F_{1,n}^2\,|\,\Cal I_1) - E{f_n}^2| + | E f^2 - E{f_n}^2|,
  \endmultline 
$$
to get (3.2) it  remains to prove
$$
   \| E(\bar F_{1,n}^2\,|\,\Cal I_1) - E{f_n}^2\|_1 \to 0. \tag{3.7}
$$

We have
$$
  \bar F_{0,n}^2 = \frac1n \sum_{j=0}^{n-1}  f_n ^2\circ T_{0,j} + \frac2n \sum_{k=1}^{n-1} \sum_{j=0}^{k-1} 
  f_n\circ T_{0,k}  f_n \circ T_{0,j}.
$$
The sigma fields $\Cal I_1$ and  $\Cal I_2$ are independent hence by Ergodic Theorem
$$
  \frac1n \sum_{j=0}^{n-1}  E(f_n^2\circ T_{0,j}\,|\,\Cal I_1) \to E[E(f_n^2\,|\,\Cal I_1)\,|\,\Cal I_2] = E[f_n^2] \tag{3.8}
$$
(in $L^1$ and a.s.).

(3.7) is thus follows from
$$
   E\Big(\frac1n \sum_{k=1}^{n-1} \sum_{j=0}^{k-1} f_n \circ T_{0,k}f_n \circ T_{0,j} \Big| \Cal I_1\Big) \to 0
  \quad \text{in}\quad L^1.   \tag{3.9}
$$
We will prove (3.9).

\bigskip

From Lemma 4 we deduce that
for the sequence of $C_n$ defined above there exist integers $\ell_k\nearrow \infty$ and real numbers 
$\epsilon_k\searrow 0$ such that for $n\to \infty$
$$\gathered
  \|f-f_n\|_2 \to 0, \quad \| E( f_n | \Cal F_{\infty,-\ell_n} \vee \Cal I_1) \|_2 <\epsilon_n, \\
  \frac{\ell_n^2C_n^4}n \to 0 \,\,\, \text{hence}\,\,\, \frac{C_n^2}n \to 0, \,\,\, \text{hence}\,\,\, C_n = o(\sqrt n).
  \endgathered \tag{3.10}
$$
To see it, notice that similarly as the $C_n$, we can let the $\ell_n$ grow arbitrarily slowly. The $\epsilon_n$ will of course decay 
slowly as well, we only need $n$ big for $\epsilon, \ell$ given. \newline
We will use this to get (3.13).
\bigskip

By Lemma 1 (see also \cite{GLV}) we can see that 

$\big(E(f_n \circ T_{0,k} \sum_{j=0}^{(k-1)\wedge (k-\ell)}  f_n \circ T_{0,j} \,|\,\Cal I_1)\big)_k$ \newline
are martingale differences  hence are mutually orthogonal. 
\bigskip

Without loss of generality we can suppose $\ell_n\geq 2$. We have
$$\gathered
  J_n = \frac1{\sqrt n} \sum_{k=1}^{n-1} \frac1{\sqrt n} f_n \circ T_{0,k} \sum_{j=0}^{k-1} f_n \circ T_{0,j} = \\
  \frac1{\sqrt n} \sum_{k=\ell_n}^{n-1} \frac1{\sqrt n}  f_n \circ T_{0,k} 
  \sum_{j=0}^{k-\ell_n}  f_n \circ T_{0,j} + \\
  \frac1{\sqrt n} \sum_{k=1}^{n-1} \frac1{\sqrt n}  f_n \circ T_{0,k} \sum_{j=(k-\ell_n+1)\vee 0}^{k-1} f_n \circ T_{0,j} = \\
  I_n + II_n. 
  \endgathered \tag{3.11}
$$
The proof of (3.9) thus reduces to 

\proclaim{Lemma 5}
$$
  E(I_n | \Cal I_1) \to 0, \quad E(II_n | \Cal I_1) \to 0 \quad \text{in}\quad L^1. \tag{3.12}
$$
\endproclaim

First, we prove the convergence for $II_n$. \newline
Because $|f_n| \leq 2C_n$ and $f_n \circ T_{0,k} 
\sum_{j=(k-\ell+1)\vee 0}^{k-1} f_n \circ T_{0,j}$ are 
martingale differences, we have
$$\gathered
  \| II_n\|_2^2   = \|\frac1n \sum_{k=1}^{n-1}  f_n \circ T_{0,k} 
  \sum_{j=(k-\ell_n+1)\vee 0}^{k-1} f_n \circ T_{0,j}\|_2^2   = \\
  \frac1{n^2} \sum_{k=1}^{n-1} \| f_n \circ T_{0,k} \sum_{j=(k-\ell_n+1)\vee 0}^{k-1} f_n \circ T_{0,j}\|_2^2 \leq \\
  \frac1{n^2} \sum_{k=1}^{n-1} \Big(\sum_{j=(k-\ell_n+1)\vee 0}^{k-1}  \| f_n \circ T_{0,k}f_n \circ T_{0,j}\|_2\Big)^2 
  \leq   16 \frac{\ell_n^2 n}{n^2} C_n^4 = 16 \frac{\ell_n^2 }{n} C_n^4 \to 0 
  \endgathered \tag{3.13}
$$
(at the end we used (3.10)).

\bigskip

Next we prove that
$$
  E(I_n | \Cal I_1)  \to 0  \quad \text{in} \quad L^1. \tag{3.14}
$$
Denote $\Cal I_{1,k} = \Cal I \cap \Cal F_{\infty, k}$.
By Lemma 1 $E( f_n \circ T_{0,k}  f_n \circ T_{0,j} \,|\,\Cal I_1)$ is $\Cal F_{\infty, k}$ measurable hence
$E( f_n \circ T_{0,k}  f_n \circ T_{0,j} \,|\,\Cal I_1) = E( f_n \circ T_{0,k}  f_n \circ T_{0,j} \,|\,\Cal I_{1,k})$.
Without loss of generality we can suppose $\ell_n\geq 2$. 
We thus have
$$\gather
  E(I_n | \Cal I_1) = \frac1n \sum_{k=\ell_n}^{n-1} \sum_{j=0}^{k-\ell_n}
  E( f_n \circ T_{0,k}  f_n \circ T_{0,j} \,|\,\Cal I_1) = \\
  E\Big[ \sum_{k=\ell_n}^{n-1} \frac1{\sqrt n}  
  E(f_n \circ T_{0,k}\,|\,\Cal I_{1,k}\vee \Cal F_{\infty, k+1-\ell_n}) 
  \frac1{\sqrt n} \sum_{j= 0}^{k-\ell_n} f_n \circ T_{0,j} \,\Big|\, \Cal I_{1,k}\Big] = \\
  \sum_{k=\ell_n}^{n-1} E(X_{n,k}  \,|\, \Cal I_{1,k})
  \endgather
$$
where
$$
  X_{n,k} = \Big(\frac1{\sqrt n} E(f_n \circ T_{0,k}\,|\,\Cal I_{1,k}\vee \Cal F_{\infty, k+1-\ell_n}) \Big)
  \Big( \frac1{\sqrt n} \sum_{j= 0}^{k-\ell_n} f_n \circ T_{0,j}\Big) = X'_{n,k}X_{n,k}''
$$ 
($\ell_n\geq 2$, $k\geq \ell_n$). \newline
Notice that for $n$ fixed, $X_{n,k}$ are $L^2$ martingale differences. \newline
We will show that $I_n$ are uniformly integrable and then using McLeish's CLT (Theorem E) we'll show that 
$I_n\to 0$ in probability. This will imply (3.14).

First we show that $J_n$ are uniformly integrable:

The $L^1$ norms of 
$$
  \frac1n \sum_{j=0}^{n-1} \sum_{k=0}^{n-1} f_n\circ T_{0,j} f_n\circ T_{0,k}
  =\Big( \frac1{\sqrt n} \sum_{j=0}^{n-1} f_n\circ T_{0,j} \Big) 
  \Big( \frac1{\sqrt n} \sum_{k=0}^{n-1} f_n\circ T_{0,k} \Big) \tag{3.15}
$$
equal $\|f_n^2\|_1 \to \|f^2\|_1$.

The random variables $f_n \circ T_{0,j}$ are martingale differences and $\|f - f_n\|_2 \to 0$ hence using Theorem F
we can deduce that $(1/\sqrt n) \sum_{j=0}^{n-1} f_n\circ T_{0,j}$ converge in law to a distribution with mean 
$\|f^2\|_1$ and characteristic function
$E \exp(\frac12 \eta^2t^2)$ where $\eta^2 = E(f^2 | \Cal I_2)$. 
Therefore the sums (3.15) are uniformly integrable.

Because $|f_n| \leq 2|f|$, from Ergodic Theorem for $(1/n) \sum_{j=0}^{n-1} (f\circ T_{0,j})^2$ we deduce uniform 
integrability of $(1/n) \sum_{j=0}^{n-1} (f_n\circ T_{0,j})^2$.

From this we deduce that \newline
$(1/n) \sum_{k=1}^{n-1} \sum_{j=0}^{k-1} f_n\circ T_{0,j} f_n\circ T_{0,k} +
(1/n) \sum_{k=0}^{n-2} \sum_{j=k+1}^{n-1} f_n\circ T_{0,j} f_n\circ T_{0,k}$ are uniformly integrable.
Because
$$
  \frac1n \sum_{k=1}^{n-1} \sum_{j=0}^{k-1} f_n\circ T_{0,j} f_n\circ T_{0,k} =
  \frac1n \sum_{k=0}^{n-2} \sum_{j=k+1}^{n-1} f_n\circ T_{0,j} f_n\circ T_{0,k},
$$ 
we deduce that 
$J_n = (1/n) \sum_{k=1}^{n-1} \sum_{j=0}^{k-1} f_n\circ T_{0,j} f_n\circ T_{0,k}$ are uniformly integrable. From this 
and uniform integrability of $II_n$ (see $(3.13)$) we deduce uniform integrability of $I_n$.
\bigskip

Next we will show that the McLeish's central limit theorem (Theorem E) applies. 
We thus have to verify that
\roster
\item"(i)" $\max_{1\leq k \leq n-1} |X_{n,k}|$ are uniformly bounded in $L^2$ norm,
\item"(ii)" $\max_{1\leq k \leq n-1} |X_{n,k}| \to 0$ in probability,
\item"(iii)" $\sum_{k=1}^{n-1} X_{n,k}^2$ converge in probability to a constant $\sigma^2=0$.
\endroster
\medskip

By (3.10) we can suppose $C_n = o(\sqrt n)$. For $n$ fixed the random variables \newline
$(1/\sqrt n) 
E(f_n \circ T_{0,k}\,|\,\Cal I_{1,k}\vee \Cal F_{\infty, k-\ell_n})$ 
are thus uniformly bounded and for $n\to\infty$ they (uniformly) converge to zero in $L^\infty$.

$X_{n,k}''$ are (for $n$ fixed) stationary martingale differences hence the weak invariance principle applies 
(cf\. \cite{HaHe}); 
 $X_n'' = \max_{k\leq n} (X_{n,k}'')^2$ are thus uniformly integrable. \newline
As shown above, $\max_k |X'_{n,k}| \to 0$ in $L^\infty$; from this and uniform integrability of $\max_{k\leq n} (X_{n,k}'')^2$
we get (i) and (ii).

 We have
$$\gather
  \sum_{j=1}^{n-1} X_{n,j}^2 =
  \frac 1n \sum_{j=1}^{n-1} (E(f_n \circ T_{0,j}\,|\, \Cal F_{\infty, j-\ell_n}\vee \Cal I_{1,k}))^2 (X_{n,j}'')^2 \leq \\
   X_n''  \frac 1n \sum_{j=1}^{n-1} E\Big((E(f_n\circ T_{0,j} \,|\, \Cal F_{\infty, -\ell_k}\vee \Cal I_{1,k}))^2.
  \endgather
$$  
By (3.10) 
$\|E(f_n\,|\, \Cal F_{\infty, -\ell_n}\vee \Cal I_1)\|_2 < \epsilon_n \to 0$ and (iii) (with $\sigma^2=0$) follows.

By Theorem E (McLeish's CLT) we thus have $I_n \to 0$ in distribution hence in probability. Because $I_n$ are uniformly 
integrable, (3.14) follows. Together with (3.13) 
this proves (3.12) hence (3.9) and thus ends the proof of Proposition 3a.

\enddemo
\qed

\demo{Proof of Proposition 3}
Proposition 3a implies assumptions of Proposition 4 hence
 the sums $(1/\sqrt{n_1n_2}) \sum_{i_1=1}^{n_1} \sum_{i_2=1}^{n_2} U_{i_1,i_2} f$ converge in law to
$\Cal N(0, \sigma^2)$ where $\sigma^2 = \|f\|_2^2$.

From this we get Proposition 3 by induction in the same way as in \cite{V15}. 
For reader's convenience we give a sketch of a proof here.

Suppose that $(1/{\sqrt{{|\underline{n_d}|}}}) \sum_{i_2=1}^{n_2} \dots \sum_{i_d=1}^{n_d} U_{0,\dots,i_d} f \to
\Cal N(0, \sigma^2)$ in law, $2\leq d$, $\min(n_2, \dots,n_d) \to \infty$, for all $f$ satisfying (1.3);
$|\underline{n_d}| = n_2\cdots n_d$.
We note that $a_0f+ a_1U_{1,0,\dots,0}f + a_mU_{m,0,\dots,0}f$ satisfies (1.3) as well hence the vector of
$F_{0,\underline{v}}, \dots, F_{m,\underline{v}}$ where $\underline{v} = (n_2, \dots,n_d)$, for
$F_{i,\underline{v}} = (1/{\sqrt{{|\underline{n_d}|}}}) \sum_{i_2=1}^{n_2} \dots \sum_{i_d=1}^{n_d} U_{i,\dots,i_d} f$ 
the vector $(F_{i,\underline{v}})_{i=0}^m$
is in distribution close to a Gaussian vector of mutually orthogonal random variables.

Therefore, $(1/m) \sum_{i=1}^m F_{i,\underline{v}}^2$ is (in $L^1$) close to $\sigma^2$ and in the same way as 
in \cite{V15} we extend this to sums $(1/N) \sum_{i=1}^N F_{i,\underline{v}}^2$ for all $N>>m$. Then we get the statement of Theorem 2 using Proposition 4 with $\eta^2=\|f\|_2^2$.
\enddemo
\qed

\demo{Proof of Theorem 2} 

1. If all transformations $T_{e_q}$ are of positive entropy in $\Cal I$ then as shown in Remark 1 after Theorem 2
the factor $\Cal I$ is the same 
case as in the Wang-Woodroofe example and we get a convergence to a non normal law.

2. Let one of the transformations $T_{e_q}$ be of zero entropy in $\Cal I$. Then by \cite{V87}
$E(f\,|\, \Cal F_{-1}^{(q)}\vee \Cal I) = 0$ and Proposition 3 applies (it is the same proof as in \cite{GLV}).
\enddemo
\qed

\comment
\underbar{Remark} In the Wang-Woodroofe example, it is in fact supposed that $T_{e_q}$ be of zero entropy in 
$\bar \Cal I_q$. We have $\Cal I = \Cal I_q \vee \bar \Cal I_q$ where $\Cal I_q$ and $\bar \Cal I_q$ are independent
$T_{e_q}$-invariant sigma algebras hence the entropy of $T_{e_q}$ in $\Cal I$ is the same as in $\bar \Cal I_q$.
\endcomment

\enddocument